\documentclass{article}
\usepackage{amsmath}
\usepackage{amsfonts}
\usepackage{amssymb}
\usepackage{graphics}
\usepackage{epsfig}
\usepackage{textcomp}

\newtheorem{lemma}{Lemma}[subsection]

\newtheorem{definition}{Definition}[subsection]
\newtheorem{notation}{Notation}[subsection]
\newtheorem{remark}{Remark}[subsection]
\title{Systole and $\lambda_{2g-2}$ of closed hyperbolic surfaces of genus $g$}
\author{Sugata Mondal}
\begin{document}

\maketitle
\begin{abstract} We apply topological methods to study eigenvalues of the  Laplacian on closed hyperbolic surfaces. For any closed hyperbolic surface $S$ of genus $g$, we get a geometric lower bound on ${\lambda_{2g-2}(S)}:$ ${\lambda_{2g-2}(S)} > {1/4} + {\epsilon_0}(S),$ where ${\epsilon_0}(S) >0$ is an explicit constant which depends only on the {\it systole} of $S$.
\end{abstract}

\textbf{Introduction.}
Here a {\it hyperbolic surface} is a complete two dimensional Riemannian manifold with sectional curvature equal to $-1$. Any hyperbolic surface is isometric to the quotient $\mathbb H/ \Gamma $, where $\mathbb H$ is the Poincar\'e upper halfplane and $\Gamma$ is a {\it Fuchsian group}, i.e. a discrete torsion-free subgroup of PSL(2,$\mathbb R$). The {\it Laplacian} on $\mathbb H$ is the differential operator which associates to a $C^2$- function $f$ the function $$ \Delta f(z) = {y^2}(\frac{{\partial^2}f}{{\partial x}^2} + \frac{{\partial^2}f}{{\partial y}^2}).$$ It induces a differential operator on $S = {\mathbb H/ \Gamma}$ which extends to a self-adjoint operator $\triangle_S$ densely defined on ${L^2}(S).$ Its domain consists of all distributions $\phi$ such that $\Delta \phi \in {L^2}(S)$. The Laplacian is a non-positive operator whose spectrum is contained in a smallest interval $( -\infty, -{\lambda_0}(S) ] \subset {\mathbb R^{-}} \cup \{0\}$ with ${\lambda_0}(S) \geq 0$. The {\it Rayleigh quotient}s allow us to characterize the {\it bottom of the spectrum} of $S,$ ${\lambda_0}(S)$: $$ {\lambda_0}(S) = inf \frac{{\int_S}{\lVert \nabla \phi \rVert}^2 dv}{{\int_S} \phi^2 dv},$$ where the infimum is taken over all non-constant smooth functions $\phi$ with compact support. Recall that the bottom of the spectrum on $\mathbb H$ is ${\lambda_0}(\mathbb H) = {1/4}$ (cf. \cite[p. 46, Theorem 5]{Cha}).

\begin{definition}
Let $\lambda > 0.$ A function $ f : S \rightarrow \mathbb R$ is a $\lambda$-{\it eigenfunction} if $ f \in {L^2}(S)$ and satisfies $ \Delta f + \lambda f = 0 .$ One calls $\lambda$ an {\it eigenvalue}. When $0 < \lambda \leq {1/4}$, $\lambda$ is called a {\it small eigenvalue} and $f$ is called a {\it small eigenfunction}.
\end{definition}
When $S$ is a compact hyperbolic surface, the {\it spectrum} of $S$ is a discrete set:
$$0={\lambda_0}(S) < {\lambda_1}(S) \leq {\lambda_2}(S)...\leq {\lambda_n}(S) \leq...$$
where in the above sequence each number is repeated according to its multiplicity as an eigenvalue and ${\lambda_i}(S)$ denotes the $i$-th non-zero eigenvalue of $S$ for $i \geq 1$.

\begin{definition} For a hyperbolic surface $S$ the {\it systole} $s(S)$ of $S$ is defined to be the minimum of the lengths of closed geodesics on $S$.
\end{definition}
The main result of this paper is the following.

\textbf{Theorem 1} {\it Let $S$ be a closed hyperbolic surface of genus $g$. Then there exists an explicit constant ${\epsilon_0}(S) > 0$, depending only on the systole of $S$ such that ${\lambda_{2g-2}}(S) > {1/4} + {\epsilon_0}(S).$}

Let $\mathcal{M}_g$ be the moduli space of closed hyperbolic surfaces of genus $g$. We recall some general facts about the behavior of $\lambda_{2g-2}$ as a function on $\mathcal{M}_g$. Any eigenvalue $\lambda_i$, in particular ${\lambda_{2g-2}}$, is a continuous function on ${\mathcal M}_g$ (see for instance \cite{C-C}). The moduli space ${\mathcal M}_g$ is the space of all closed hyperbolic surfaces of genus $g$ up to isometry. Recall that the set ${{\mathcal I}_\epsilon}= \{ S \in {{\mathcal M}_g}: s(S) \geq \epsilon \}$ is compact (\cite[p. 163]{Bu}). By \cite{O-R} ${\lambda_{2g-2}}(S) > {1/4}$ for all $S \in {\mathcal M}_g$. Hence there exists a non-zero constant ${\eta}(\epsilon)$ such that ${\lambda_{2g-2}}(S) > {1/4} + {\eta}(\epsilon)$ for all $S \in {{\mathcal I}_\epsilon}$. This proves the Theorem with ${\epsilon_0}(S)= {\eta}(s(S)).$ The content of Theorem 1 is to make this constant explicit in terms of the geometry of $S$. We shall see that we can take ${\epsilon_0}(S)$ to be any positive number smaller than $$\textrm{min}\displaystyle \{ \frac{1}{4(g-1)}, {\frac{1}{4}}({({\frac{\cosh \rho_0}{\sinh \rho_0}})^2} -1) \}$$ where $2 s(S) {\sinh \rho_0} = |S|$.

We now briefly sketch the proof of the above theorem. It uses topological methods as in \cite{O-R}. First we recall that an open subset of a surface $S$ is called {\it incompressible} if the fundamental group of any of its connected components maps injectively into ${\pi_1}(S).$ We use the convention that simply connected open subsets of $S$ are incompressible. Let ${\mathcal E}_{\lambda}$ denote the eigenspace of the Laplacian on $S$ for the eigenvalue $\lambda.$ For $\epsilon > 0$, let ${\mathcal E}^{\frac{1}{4} + \epsilon}$ be the direct sum of eigenspaces ${\mathcal E}_{\lambda}$ with $\lambda \leq {\frac{1}{4} + \epsilon}.$ For $f\neq 0$ $\in {\mathcal E}^{\frac{1}{4} + \epsilon}$, define the {\it nodal set} ${\mathcal Z}(f)$ as ${{f}^{-1}}(0).$ It is proved in \cite{O-R}, using the analyticity of eigenfunctions on $\mathbb H$, that ${\mathcal Z}(f)$ is the union of a finite graph and a discrete set. Let ${\mathcal G}(f)$ be the subgraph of ${\mathcal Z}(f)$ obtained by suppressing those connected components which are homotopic to a point on $S$ (equivalently, those which are contained in a topological disc). Due to this modification, each component of $S \setminus {\mathcal G}(f)$ is incompressible.  One of the main observation in \cite{O-R} was that for any $f\neq 0$ $\in {\mathcal E}^{\frac{1}{4}}$, the Euler characteristic of at least one component of $S \setminus {\mathcal G}(f)$ is strictly negative. For $\epsilon > 0$ there is no reason, in general, to believe such a result for $f\neq 0$ $\in {\mathcal E}^{\frac{1}{4}+\epsilon}.$ However, we will show the following

\textbf{Lemma 1} {\it Let $S$ be a closed hyperbolic surface of genus $g$. Then there exists an explicit constant ${\epsilon_0}(S) > 0$ depending only on the genus $g$ and the systole of $S$, such that for any $f\neq 0$ $\in {\mathcal E}^{\frac{1}{4} + {\epsilon_0}(S)}$, the Euler characteristic of at least one component of $S \setminus {\mathcal G}(f)$ is strictly negative.}

Let $\Sigma$ be a Riemannian surface. Let $\varOmega \subseteq \Sigma$ be an open set such that the closure $\bar{\varOmega}$ is a submanifold with piecewise smooth boundary. Then denote by $\Delta$ the Laplace operator of $\Sigma$ restricted to $\varOmega$. {\it Dirichlet eigenvalue}s of $\varOmega$ are the $\lambda$'s such that the problem:
\begin{equation*}
\left.
\begin{aligned}
\Delta u & = {\lambda u} ~& \mbox{ on } & ~\varOmega, \\
u & = 0 ~& \mbox{ on } & ~\partial \varOmega.
\end{aligned}
~~~ \right\} ~~~~~
\end{equation*}\\
admits a non-zero solution $u$, continuous on $\bar{\varOmega}$ and smooth on $\varOmega$. The smallest $\lambda$ for which such a solution exists is denoted by $\lambda_0 (\varOmega)$ and is called {\it the first Dirichlet eigenvalue} of $\varOmega$. This number can be defined in terms of Rayleigh quotients in a similar way as the bottom of the spectrum of $\varOmega$, $$ {\lambda_0}(\varOmega) = inf \frac{{\int_{\varOmega}}{\lVert \nabla \phi \rVert}^2 dv}{{\int_{\varOmega}} \phi^2 dv},$$ where the infimum is taken over all non-zero smooth functions $\phi$ with compact support in $\varOmega$. From this characterization it is evident that for any two submanifolds $\varOmega_1$ and $\varOmega_2$ as above with compact closure, then ${\lambda_0}(\varOmega_1) \gneq {\lambda_0}(\varOmega_2)$ when $\varOmega_1 \subsetneq \varOmega_2$. The above Lemma will be deduced from the following:

\textbf{Proposition 1} {\it Let $S$ be a closed hyperbolic surface of genus $g$. Let $\Omega \subset S$ be a surface with smooth boundary which is homeomorphic either to a disc or to an annulus. Then there exists a constant ${\epsilon}(\Omega) > 0$ depending on the length $l_{\Omega}$ of the geodesic in $S$ homotopic to a generator of ${\pi_1}(\Omega)$ and the area of $\Omega$ such that the first Dirichlet eigenvalue of $\Omega$ satisfies: ${\lambda_0}(\Omega) > {1/4} + {\epsilon}(\Omega).$ Furthermore there exists an explicit constant ${\epsilon_0}(S) > 0$ depending only on the systole of $S$ such that ${\epsilon}(\Omega) > {\epsilon_0}(S).$}

\begin{notation} For any surface $\Omega \subseteq S $ with smooth boundary, $|\Omega|$ denotes the area of $\Omega$ for the area measure on $S$ and $L(\partial \Omega)$ denote the length of the boundary of $\Omega$.
\end{notation}
We shall see in the proof that $\epsilon (\Omega)$ is a strictly decreasing function of $|\Omega|$ when  $l_{\Omega}$ is kept fixed and a strictly increasing function of $l_{\Omega}$ when $|\Omega|$ is kept fixed. The statement in the proposition then follows from the observation that both the parameters i.e. $|\Omega|$ and  $l_{\Omega}$ are bounded: the first one being bounded above by $4\pi(g-1)$ and the last one being bounded below by $s(S)$.

The proof of the above proposition depends mainly on two geometric inequalities: the Faber-Krahn isoperimetric inequality and the Cheeger's inequality. The scheme of the proof of Theorem 1 follows then the same lines as the one of Theorem 1 in \cite{O-R}.

Existence of surfaces with small eigenvalues was proved originally by B. Randol \cite{R1} using the famous trace formula of A. Selberg. We shall recall another method of P. Buser \cite{Bu} for the construction of such surfaces using max-min principle. The construction is carried out by first considering a genus $g$ hyperbolic surface admitting a pair of pants decomposition with very short boundary geodesics, then constructing an orthogonal family of functions with small Rayleigh quotient. The number of functions in that family is exactly $(2g-2)$. This gives the existence of at least $(2g-3)$ small eigenvalues (which is the maximum possible number by \cite{O-R}).

After proving theorem 1 in the first section, in the second part of the paper we study the behavior of $\lambda_i$ as a function on the {\it moduli space} ${\mathcal M}_g$. We recall that the moduli space ${\mathcal M}_g$ is the space of all closed hyperbolic surfaces of genus $g$ up to isometry. We focus our interest on the first $2g-2$ non-zero eigenvalues. Theorem 1 (or even a continuity argument on ${\mathcal M}_g$) implies one direction of the following

\textbf{Claim 1} {\it For a family $S_n$ of compact hyperbolic surfaces in ${\mathcal M}_g$, ${\lambda_{2g-2}}(S_n)$ tends to ${1/4}$ if and only if the systole $s(S_n)$ tends to zero.}

The other direction follows from a construction due to P. Buser \cite{Bu}.

The above proposition can be compared with the following result of Schoen, Wolpert and Yau \cite{S-W-Y}. Let $M$ be a closed oriented surface of genus $g$ with a metric of (possibly variable) Gaussian curvature $K$. For an integer $n \geq 1$ consider the family ${\tilde C}_n$ of curves on $M$ which are disjoint union of simple closed geodesics and which divide $M$ into $n+1$ components (necessarily $n \leq 2g-3$). Define a number $l_n$ by
$$ l_n = min \{ L(C) : C \in {\tilde C}_n \}.$$ where $L(C)$ denotes the length of $C$. Then

\textbf{Theorem (Schoen-Wolpert-Yau).} {\it Suppose for some constant $k > 0$ we have $-1 \leq K \leq -k.$ Then there exist positive constants $\alpha_1, \alpha _2 $ depending only on $g$ such that for $1\leq n \leq 2g-3$, we have ${\alpha_1} {{k}^{3/2}} l_n \leq \lambda_n \leq {\alpha_2}{l_n}$ and ${\alpha_1}k \leq \lambda_{2g-2} \leq \alpha_2.$}

Recall that the Bers constant $\beta$ \cite{B} which depends only on $g$ has the property that ${l_{2g-2}} < \beta$. So this theorem implies that ${\lambda_{2g-2}}$ is bounded above by a constant depending only on $g$. Observe also that the Buser's construction (\cite[Theorem 8.1.3]{Bu}) leads to the same conclusion. Namely by Buser's construction for any $\delta >0$ there exists a constant $\epsilon > 0$ such that  ${\lambda_{2g-2}} < \frac{1}{4} + \delta$ for any $S \in {\mathcal M}_g$ with $s(S) < \epsilon.$ Since $\lambda_{2g-2}$ is a continuous function on ${\mathcal M}_g$ and ${{\mathcal I}_\epsilon}= \{ S \in {{\mathcal M}_g}: s(S) \geq \epsilon \}$ is compact the existence of an upper bound is clear. In this context we would like to mention a paper due to Dodziuk, Pignataro, Randol and Sullivan \cite{D-P-R-S} where the authors obtained a result similar to the one of \cite{S-W-Y} in the context of possibly non-compact hyperbolic surfaces.

In \S 2 we will study the behavior of ${\lambda_i}(S)$ as $s(S)$ tends to zero. More precisely, let $\overline{{\mathcal M}_g}$ denote the compactification of ${\mathcal M}_g$ obtained by adding the moduli spaces of (not necessarily connected) non-compact finite area hyperbolic surfaces with area equal to $4 \pi(g-1)$. Let ${\partial {\mathcal M}_g}= \overline{{\mathcal M}_g} \setminus {\mathcal M}_g$ be the corresponding boundary of ${{\mathcal M}_g}$. We study the behavior of ${\lambda_i}(S_n)$ when ${S_n} \in {\mathcal M}_g$ tends to a point in ${\partial {\mathcal M}_g}$. By the above theorem of Schoen, Wolpert and Yau and the discussion after, it is clear that ${\lambda_i}(S)$ is bounded above for all $S \in {\mathcal M}_g$ and for $1 \leq i \leq 2g-2$. Indeed the method using Buser's construction works for any $i$, showing that $\lambda_i$ is bounded by a constant depending only on $g$ and $i$. So for any $i$ we can consider the set $$V_i = \{ {\lim_{n \rightarrow \infty}} {{\lambda_i}(S_{n})} :({S_n}) ~ \textrm{is a sequence in} ~ {{\mathcal M}_g} ~ \textrm{converging to a point in} ~ {\partial {\mathcal M}_g} $$ $$ \textrm{such that} {\lim_{n \rightarrow \infty}} {{\lambda_i}(S_{n})}  ~ \textrm{ exists}\}.$$ With this notation, the above claim says that ${V_{2g-2}} = \{ \displaystyle \frac{1}{4} \}.$ We next prove

\textbf{Claim 2} {\it For any $1 \leq i \leq 2g-3$, there exists a ${\Lambda_i}(g) \in ( 0, \displaystyle \frac{1}{4}]$ such that $V_i$ contains the interval $[0, {\Lambda_i}(g)]$.}

We shall use a result of Courtois and Colbois \cite[Theorem 0.1]{C-C} to prove this claim.

In \S 3 we study non-compact hyperbolic surfaces of finite area. Recall that for a non-compact  hyperbolic surface $S$ of finite area, the spectrum of the Laplace operator is composed of two parts: the discrete part and the continuous part. The continuous part covers the interval $[{1/4}, \infty)$ and is spanned by Eisenstein series. The discrete part is the union of the {\it residual spectrum} and {\it cuspidal spectrum}. The residual spectrum is a finite set contained in the interval $(0,{1/4})$ and it corresponds to poles of the analytic continuation of Eisenstein series. The cuspidal spectrum consists of those eigenvalues whose associated eigenfunctions tend to zero uniformly near any cusp. The number of cuspidal eigenvalues is known to be infinite for arithmetic groups \cite{I}. The cuspidal eigenvalues can possibly appear anywhere in the interval $(0, \infty).$ Denote ${{\lambda^c}_i}(S)$ the $i$-th cuspidal eigenvalue of $S$.

In analogy to Theorem 1, one can investigate

\textbf{Conjecture.} {\it Let $S$ be a finite area hyperbolic surface of type $(g, n)$. Then there exists an explicit constant ${\epsilon_0}(S) > 0$, depending only on the systole of $S,$ such that ${{\lambda^c}_{2g-2+n}}(S) > 1/4 + {\epsilon_0}(S)$.}

This would be an extension of a result of Jean-Pierre Otal and Eulalio Rosas (Theorem 2 in \cite{O-R}). However our methods do not suffice to settle this conjecture. In this connection we state the following conjecture of Jean-Pierre Otal and Eulalio Rosas in \cite{O-R} which is motivated by \cite[Prop. 2, Prop. 3]{O}

\textbf{Conjecture.} {\it Let $S$ be a finite area hyperbolic surface of type $(g, n)$. Then ${{\lambda^c}_{2g-2}}(S) > 1/4$.}

Now we consider a finite area hyperbolic surface $S$ of type $(g,n)$. Denote by ${\mathcal T}_{g,n}$ the {\it Teichm\"{u}ller space} of all marked hyperbolic surfaces of type $(g,n)$. For any choice of pair of pants decomposition of $S$ one can define a system of coordinates on ${\mathcal T}_{g,n}$, the {\it Fenchel-Nielsen} coordinates which consists, for each curve in the pants decomposition, of the length of that curve and a {\it twist} parameter along that curve (\cite[Chapter 6]{Bu}). Now we consider the set ${{{\mathcal T}^0}_{g,n}}$ of all hyperbolic surfaces in ${\mathcal T}_{g,n}$ for which all twist parameters are equal to zero. Each surface in ${{{\mathcal T}^0}_{g,n}}$ carries an involution $\iota$ which when restricted to each pair of pants is the orientation reversing involution that fixes the boundary components. This involution induces an involution on each eigenspace of the Laplacian. The eigenfunctions corresponding to the eigenvalue $-1$ are called {\it antisymmetric} and the corresponding eigenvalue is called an {\it antisymmetric eigenvalue}. We denote the i-th antisymmetric cuspidal eigenvalue of $S$ by ${{{\lambda}^{o,c}}_{i}}(S).$

\textbf{Theorem 2} {\it For every surface $S \in {{{\mathcal T}^0}_{g,n}}$ there exists an explicit constant ${\epsilon_0}(S) > 0,$ depending only on the systole of the surface $S,$ such that ${{{\lambda}^{o,c}}_{g}}(S) > {1/4} + {\epsilon_0}(S).$}

Indeed, the constant ${\epsilon_0}(S)$ can be taken equal to any number below $$\textrm{min}\displaystyle \{ \frac{1}{2(2g-2+n)}, {\frac{1}{4}}({({\frac{\cosh \rho_0}{\sinh \rho_0}})^2} -1) \}$$ where $2 s(S) {\sinh \rho_0} = |S|.$

\subsubsection*{\centerline{Acknowledgement}} The author would like to express his sincere gratitude to his advisor Jean-Pierre Otal whose encouragement, kindness and patience were significant ingredients in the work. The author was supported during this research by the Indo-French CEFIPRA-IFCPAR grant.
\section{Geometric lower bound on ${\lambda_{2g-2}}(S)$} In this section we shall prove Theorem 1. We begin by Proposition 1, which we now recall

\textbf{Proposition 1} {\it Let $S$ be a closed hyperbolic surface of genus $g$. Let $\Omega \subset S$ be a surface with smooth boundary which is homeomorphic either to a disc or to an annulus. Then there exists a constant ${\epsilon}(\Omega) > 0$ depending  on the length $l_{\Omega}$ of the geodesic in $S$ homotopic to a generator of ${\pi_1}(\Omega)$ and the area of $\Omega$ in $S$ such that the first Dirichlet eigenvalue of $\Omega$ satisfies: ${\lambda_0}(\Omega) > {1/4} + {\epsilon}(\Omega).$ Furthermore there exists an explicit constant ${\epsilon_0}(S) > 0$ depending only on the systole of $S$ such that ${\epsilon}(\Omega) > {\epsilon_0}(S).$}

\textbf{ Proof.} Suppose first that $\Omega \subseteq S$ is a disc or more generally a domain such that ${\pi_1}(\Omega)$ maps to zero in ${\pi_1}(S)$. Then choose an isometric lift of $\Omega$ to $\mathbb H$, still denoted by $\Omega$. We will use the Faber-Krahn inequality (\cite[p. 87]{Cha}) in the following form:

\textbf{Theorem (Faber-Krahn inequality)} {\it Let $\Omega \subseteq {\mathbb H}$ be a domain such that $\partial \Omega$ is smooth. Let $D$ be a geodesic disc in $\mathbb H$ with same area as $D$, i.e. $ |\Omega| = | D|.$ Then, $$ {\lambda_0}(\Omega) \geq {\lambda_0}(D),$$ with equality if and only if $\Omega$ is isometric to $D$.}

Let $B(t)$ be the geodesic disc in ${\mathbb H}$ with radius $t$. The geodesic disc with same area as $\Omega$ has radius $\displaystyle t(\Omega) = 2 {{\sinh}^{-1}}(\frac{|\Omega|}{4\pi})$. By the Faber-Krahn inequality ${\lambda_0}(B(t(\Omega))) \leqslant {\lambda_0}(\Omega).$

Since $\Omega$ is contained in $S$ whose area equals $2\pi (2g-2)$, by Gauss-Bonnet theorem, $|\Omega| < 2\pi (2g-2)$. Therefore, $B(t(\Omega))$ is contained in the disc with radius $t_0 = 2 {{\sinh}^{-1}}(g-1)$. Recall that for two subsurfaces $D_1$ and $D_2$ in $\mathbb H$ with compact closure, ${\lambda_0}(D_1) > {\lambda_0}(D_2)$ when $D_1 \subsetneq D_2$. Thus ${\lambda_0}(B(t))$ is a strictly decreasing function of $t$. Hence ${\lambda_0}(B(t(\Omega))) > {\lambda_0}(B(t_0))$. Now by Theorem 5 in \cite{Cha}, we have $${\lambda_0}(B(t)) > \lim_{s \rightarrow \infty} {\lambda_0}(B(s)) = \frac{1}{4}.$$ Hence we finally have a strictly positive ${\epsilon_1}(|\Omega|)$ which depends only on the area $|\Omega|$ of $\Omega$ such that ${\lambda_0}(B(t(\Omega))) = \frac{1}{4} + {\epsilon_1}(|\Omega|).$ Since ${\lambda_0}(B(t))$ is a strictly decreasing function of $t$, ${\epsilon_1}(|\Omega|)$ is a strictly decreasing function of $|\Omega|$ which is bounded below by the constant ${\epsilon_1}(S) = {\lambda_0}(B(t_0)) - \displaystyle \frac{1}{4}$.

Suppose now that $\Omega$ is an annulus and that the image of ${\pi_1}(\Omega)$ in ${\pi_1}(S)$ is a non-trivial cyclic subgroup $\langle \tau \rangle$. Let $\mathbb T$ denote the {\it cylinder} $\mathbb H / \langle \tau \rangle$. Let $\gamma$ denote the core geodesic of $\mathbb T$ and $l$ the length of $\gamma$. Then $l$ is the length of the shortest geodesic of $S$ homotopic to a generator of ${\pi_1}(\Omega)$. Consider an isometric lift of the annulus $\Omega$ to $\mathbb H / \langle \tau \rangle$, still denoted by $\Omega$. We need to prove that ${\lambda_0}(\Omega) > \frac{1}{4} + {\epsilon_0}(S)$ where ${\epsilon_0}(S)$ depends only on $l$ and $|\Omega|$. We will use Cheeger's inequality (\cite[p. 95]{Cha}) in the following form:

\textbf{Theorem (Cheeger inequality)} {\it Let $\Omega \subsetneq \mathbb T$ be a submanifold with piecewise smooth boundary. Let $h(\Omega)$ be the {\it Cheeger constant} of $\Omega$. Then $${\lambda_0}(\Omega) \geqslant \frac{h^2(\Omega)}{4}.$$}

Recall that the {\it Cheeger constant} of $\Omega$ is equal to inf$\displaystyle \{ \frac{L(\partial V)}{|V|} \}$ where $V$ ranges over all compact submanifolds of $\Omega$ with smooth boundary.

The proof of Proposition 1 in the case of an annulus follows from Cheeger inequality and the next
\begin{lemma}\label{annuli}
Let $\Omega \subsetneq \mathbb T$ be a submanifold with piecewise smooth boundary and $h(\Omega)$ be the {\it Cheeger constant} of $\Omega$. Then: $$h(\Omega) > 1+ {\epsilon_2}(|\Omega|, l),$$ for some constant ${\epsilon_2}(|\Omega|, l) > 0$, depending only on the area of $\Omega$ and the length $l$ of the core geodesic of $\mathbb T$.
\end{lemma}
\textbf{Proof.} First we observe that the Cheeger constant is bounded below by the quantity $\displaystyle \textrm{inf} \{ \frac{L(\partial V)}{|V|} \}$ where $V$ ranges over connected submanifolds of $\Omega$. Secondly, this infimum is the same when $V$ ranges over all discs or {\it essential} annuli contained in $\Omega$. Recall that an annulus is essential when it is not homotopically trivial in $\mathbb T.$ This is because any connected, compact submanifold $V \subseteq \Omega$ is diffeomorphic either to a disc with some discs removed or to an essential annulus with some discs removed. In both cases, taking the union of $V$ with those removed discs, one obtains a submanifold $V^{'}$ which is either a disc or an essential annulus which satisfies: $L(\partial {V^{'}}) \leq L(\partial V)$ and $|V^{'} \geq |V|$. Therefore $\displaystyle \frac{L(\partial {V^{'}})}{|V^{'}|} < \frac{L(\partial V)}{|V|}.$

Suppose now that $V\subseteq \Omega$ is diffeomorphic to a disc. By the isoperimetric inequality (\cite[p. 11]{B-Z}), one has$$(\frac{L(\partial V)}{|V|})^2 \geq 1 + \frac{4 \pi}{|V|}.$$ Therefore if $V \subseteq \Omega$ then $\displaystyle (\frac{L(\partial V)}{|V|})^2 \geq 1 + \frac{4 \pi}{|\Omega|}.$

Since $|V| < 2\pi (2g-2)$, we get $\displaystyle (\frac{L(\partial V)}{|V|})^2 > 1 + \frac{1}{g-1}.$

Now we suppose that $V \subseteq \Omega$ is an essential annulus. In order to prove the claim in this case we will need the following notion of {\it symmetrization}, which is close to the notion of {\it Steiner symmetrization} (\cite[p. 18]{H}).
\begin{definition} Let $V \subseteq \mathbb T$ be an essential annulus. The {\it symmetrization} of $V$ is the annulus $V_0 \subseteq \mathbb T$ symmetric with respect to $\gamma$ with constant width and which has the same area $V$.
\end{definition}
Recall that the {\it Fermi coordinates} on $\mathbb T$ assign to each point $p$ the pair $(s, r) \in \{ \gamma \} \times \mathbb R$, where $r$ is the signed distance of $p$ from $\gamma$ and $s$ is the point of $\gamma$ nearest to $p$. After parametrizing the geodesic $\gamma$  by arc-length, these coordinates provide a diffeomorphism between $\mathbb T$ and ${\mathbb R / l \mathbb Z} \times \mathbb R $. The hyperbolic metric in these coordinates equals $ d{r^2} + {\cosh {r}}^2 d{s^2}$.

\centerline{\includegraphics[height=4in]{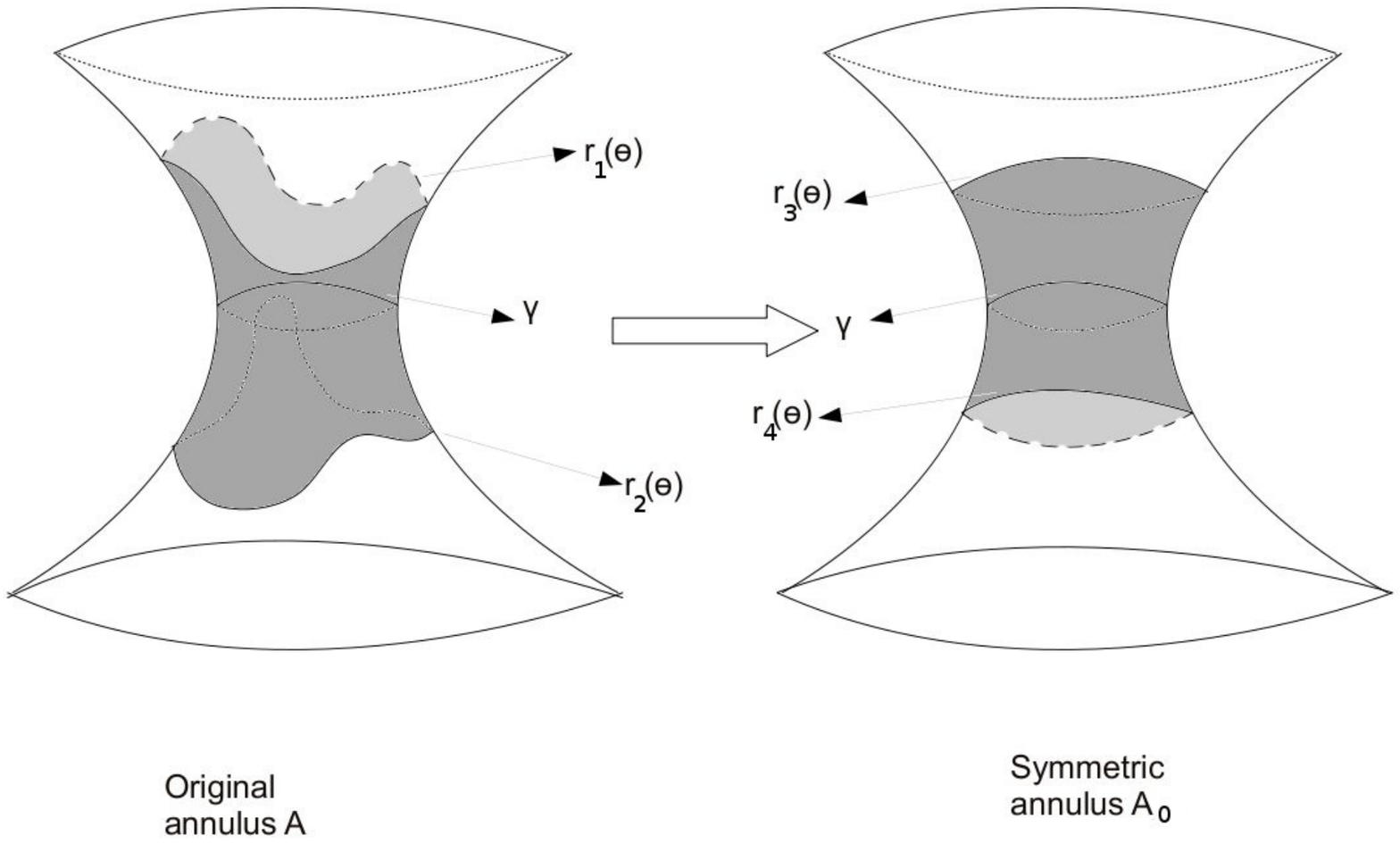}}

\begin{lemma}\label{symmetric}
Let $V \subseteq \mathbb T$ be an essential annulus with piecewise smooth boundary and $V_0$ be the symmetrization of $V$. Then $L(\partial V) \geq L(\partial V_0).$
\end{lemma}
\textbf{Proof.} First we consider the case when each component of $\partial V$ is a graph over $\gamma$. By that we mean that there exist two functions $r_1$ and $r_2$: $[0,l] \rightarrow \mathbb R$ such that $r_i$ is a piecewise smooth map (there is a partition $ 0 = {s_1} < {s_2} <  ... < {s_ m} = l$ such that each restriction $r_i | [s_j, s_{j+1}]$ is smooth) with ${r_i}(0) = {r_i}(l)$ and the components of $\partial V$ are parametrized in Fermi coordinates as: $\{ (s, {r_i}(s)), s \in [0,l] \}$ for $i = 1,2$. Then the components of the symmetrization $V_0$ of $V$ are the graphs of the constant functions $r_3 = \rho $ and $r_4 = - \rho$ with $\displaystyle \rho = {{\sinh}^{-1}}(\frac{|V|}{2l})$. Up to exchanging $r_1$ and $r_2$, we may suppose that ${r_1}(s) > {r_2}(s)$ for all $0 \leq s\leq l$. Then we calculate the areas of $V$ and $V_0$: $$|V| =  {\int_{0}^{l}} {\int_{r_1(s)}^{r_2(s)}} \cosh {r} dr ds = {\int_{0}^{l}} \{\sinh{r_2(s)} - \sinh{r_1(s)}\} ds$$
and $$|V_0| = {\int_{0}^{l}} {\int_{-\rho}^{\rho}} \cosh {r} dr ds = {\int_{0}^{l}} 2.\sinh {\rho} ds = 2l\sinh{\rho}.$$
The length of $\partial V_0$ is $$L(\partial V_0) =  2l \cosh {\rho}$$ and the length of $\partial V$ satisfies $$L(\partial V) = {\int_{0}^{l}} \{ \dot{r_1} (s)^2 + 1 \}^{1/2} \cosh {r_1(s)} ds + {\int_{0}^{l}} \{ {\dot{r_2} (s)}^2 + 1 \}^{1/2} \cosh {r_2(s)} ds $$ $$\geq {\int_{0}^{l}}\{ \cosh {r_1(s)} + \cosh{r_2(s)} \} ds.$$

Call $L_0$ the constant equal to the last expression. Observe that $L(\partial V) = L_0$ if and only if ${\dot {r_1}(\theta)}$ = 0 = ${\dot{r_2}(\theta)}$. This implies that $r_1$ and $r_2$ are constants.

One has: $${L(\partial V)}^2 -|V|^2  \geq (L_0+ |V|)(L_0 - |V|).$$
Now, $$L_0 + |V| = {\int_{0}^{l}} ((\cosh{r_2(s)} + \sinh{r_2(s)} ) + ( \cosh{r_1(s)} - \sinh{r_1(s)} )) ds$$
$$= {\int_{0}^{l}} ( \exp (r_2(s)) + \exp (- r_1(s)) ) ds$$
and similarly $$L_0 -|V| = {\int_{0}^{l}} ( \exp ( -{r_2}(s)) + \exp (r_1(s)) ) ds.$$

 Thus we have

$(L_0 + |V|)(L_0 - |V|)$
$$= \bigg({\int_{0}^{l}} ( \exp \bigg(r_2(s)\bigg) + \exp \bigg(-{r_1}(s)\bigg) ) ds\bigg) \bigg({\int_{0}^{l}} (\exp \bigg(-{r_2}(s)\bigg) + \exp \bigg({r_1}(s)\bigg) ) ds\bigg)$$

$$\geq ( {\int_{0}^{l}} (\exp \bigg({r_2}(s)\bigg) + \exp \bigg({-r_1}(s)\bigg))^{\frac{1}{2}} (\exp \bigg(-{r_2}(s)\bigg) + \exp \bigg({r_1}(s)\bigg) )^{\frac{1}{2}} ds)^2 $$
$$ \textrm{ (by H\"{o}lder's inequality)}$$

$$= ({\int_{0}^{l}} ( 2 + 2\cosh \bigg({r_1}(s) + {r_2}(s)\bigg) )^{\frac{1}{2}} ds)^2.$$
$\\$ Since $\cosh x \geq $ 1  $ \forall x$, we get  ($L_0$ + $|V|$)($L_0$ - $|V|$) $\geq 4{l^2} = {L(\partial V_0)}^2 - {A(V_0)}^2$.  Equality holds if and only if $r_1, r_2 $ are independent of $s$ and if $r_1 = - r_2$.
$\\$ Since by construction $|V|$ = $|V_0|$, the lemma is proven when $V$ is an annulus whose boundary components are graphs over $\gamma$.

Now we consider the case of an arbitrary annulus with piecewise smooth boundary. By approximation, it suffices to prove Lemma \ref{symmetric} for those $V$ which satisfy the following property: there exists a partition of $\gamma$: $0= s_1 < s_2 <... <s_k =l =0$ such that over each interval [$s_i, s_{i+1}$], $\partial V$ is the union of graphs of finitely many functions. We consider now such an annulus.  We consider the strip over [$s_i,s_{i+1}$] in $\mathbb T$ which is diffeomorphic to [$s_i,s_{i+1}$]$\times \mathbb R$ in Fermi coordinates. Denote by $V^i$ the intersection of $V$ with this strip. Let for $1 \leq i \leq k$, we denote by $ f_j, j= 0,1,2...,l(i)$ the boundary curves of $V^i$ i.e. in Fermi coordinates the components of $\partial V^{i}$ are parametrized as $\{ (s, {f_j}(s)): s \in [s_i, s_{i+1}]\}$ for $j= 0,1,2...,l(i)$ and for any $s \in [s_i, s_{i+1}]$, $r({f_0}(s)) > r({f_1}(s))>...> r({f_{l(i)}}(s))$. Now we calculate the area of $V^i$
$$|V^i| = {{\sum_{j =l(i)-1, l(i)-3,..., 1}}} {\int_{s_{i}}^{s_{i+1}}}{\int_{{f_j}(s)}^{{f_{j+1}}(s)}} \cosh r dr ds = {\sum_{j=1}^{l(i)}} {\int_{s_i}^{s_{i+1}}} {(-1)^{j+1}} \sinh f_j (s) ds.$$
The length of $\partial V^i$ is given by $$L(\partial V^i) = {\sum_{j =1}^{l(i)}} {\int_{s_{i}}^{s_{i+1}}} \{ \dot{f_j} (s)^2 + 1 \}^{1/2} \cosh {f_j(s)} ds \geq {\int_{s_{i}}^{s_{i+1}}}{\sum_{j =1}^{l(i)}} \cosh {f_j(s)}.$$
Call ${L_0}(i)$ the constant equal to the last expression and calculate
$${L(\partial V)}^2 - |V|^2 = ({\sum_{i}} L(\partial V^i))^2 - ({\sum_{i}} |V^i|)^2 \geq ({ \sum_{i}}{{L_0}(i)})^2 - ({\sum_{i}}|V^i|)^2 $$
$$=({\sum_{i}}{\int_{s_i}^{s_{i+1}}} ({\sum_{j=1}^{l(i)}} \exp[ {(-1)^{j+1}}{f_j}(s)] ds) \times({\sum_{i}} {\int_{s_i}^{s_{i+1}}} ({\sum_{j=1}^{l(i)}} \exp[ {(-1)^j}{f_j}(s)] ds)$$
$$\geq ({\sum_{i}}{\int_{s_i}^{s_{i+1}}} (\exp[ {(-1)^{0 +1}}{f_{0}}(s)] + \exp[ {(-1)^{1 +1}}{f_{1}}(s)]) ds)$$
$$\times ({\sum_{i}} {\int_{s_i}^{s_{i+1}}} (\exp[ {(-1)^{0}}{{f_{0}}(s)}] + \exp[ {(-1)^{1}}{f_{1}}(s)])ds)$$
$$ \geq  ( {\int_{0}^{l}} ( 2 + 2 \cosh ({{f_1}(s)-{f_0}}(s))^{\frac{1}{2}} ds)^2$$ $$\textrm{( by H\"{o}lder's inequality)}$$
$$\geq 4{{l}^2} = {L(\partial V_0)}^2 - |V_0|^2.$$
$\\$ Hence using the same argument as before we finally prove Lemma \ref{symmetric}.  $\Box$

So now we have $(\displaystyle \frac{L(\partial V)}{|V|}) \geq (\displaystyle \frac{L(\partial V_0)}{|V_0|}) = \displaystyle \frac{\cosh \rho}{\sinh \rho}$ where $|V| = 2l \sinh \rho$. Thus we conclude the proof of Lemma \ref{annuli} by taking $${\epsilon_2}(\Omega, l) = \displaystyle \frac{1}{2} \textrm {min} \{ \frac{\cosh \theta}{\sinh \theta} -1, {(1+ \frac{4\pi}{|\Omega|})^{\frac{1}{2}}} -1 \}$$ where $|\Omega| = 2l \sinh \theta$. $\Box$

Since $\displaystyle \frac{\cosh \rho}{\sinh \rho}$ is a strictly decreasing function of $\rho$ we have $$(\displaystyle \frac{L(\partial V)}{|V|}) \geq \displaystyle \frac{\cosh \rho_1}{\sinh \rho_1} > \displaystyle \frac{\cosh \rho_0}{\sinh \rho_0}$$ where $2l \sinh \rho_1 = |\Omega|$ and $2 s(S) \sinh \rho_0 = |S| =4\pi(g-1)$ (since $V \subseteq \Omega \subsetneq S$). To conclude the proof of Proposition 1 we take $$\epsilon_0 (S) = \displaystyle \frac{1}{2} \textrm{min}\{\epsilon_1 (S), \frac{1}{4(g-1)},  \frac{1}{4} ((\frac{\cosh \rho_0}{\sinh \rho_0})^2 -1) \}.$$$\Box$
\begin{remark}\label{rmrk}
From the expression of  $\epsilon_0(S)$ we observe that if $(S_n)$ be a sequence in ${\mathcal M}_g$, then  $\epsilon_0(S_n)$ tends to zero only if  $s(S_n)$ tends to zero. The computations in the proposition also show that for any $\Omega \subseteq S$ diffeomorphic to a disc or to an annulus one has $${\lambda_0}(\Omega) \geq \displaystyle \frac{1}{4} + 2{\epsilon_0}(S).$$
\end{remark}
\subsection{Proof of Theorem 1:}
The proof at this point follows the same lines as that of  ${\lambda_{2g-2}}(S) > \displaystyle \frac{1}{4} $ in \cite {O-R} and we refer \cite {O-R} for the details. We take ${\epsilon_0}(S)$ as in Proposition 1. Consider the space ${\mathcal E}^{\frac{1}{4} + {\epsilon_0}(S)}$. Recall that ${\mathcal E}^\lambda$ is the direct sum of the eigenspaces of the Laplacian with eigenvalues less than or equal to $\lambda$. Let $f \neq 0 \in {\mathcal E}^{\frac{1}{4} + {\epsilon_0}(S)}$. The {\it nodal set} ${\mathcal Z}(f)$ of $f$ is defined as ${{f}^{-1}}(0).$ Recall that ${\mathcal G}(f)$ is the subgraph of ${\mathcal Z}(f)$ obtained by suppressing those connected components which are zero homotopic on $S$. Each component of $S \setminus {\mathcal G}(f)$ is an open surface, may be equal to $S$ when ${\mathcal G}(f)$ is empty. The sign of $f$ on a component of $S \setminus {\mathcal G}(f)$ can be defined as follows. There is a finite collection of disjoint closed topological discs $(D_j)$ with $\partial{D_j}\bigcap {\mathcal Z}(f) = \phi$ such that each component of ${\mathcal Z}(f)$ which is zero homotopic is contained in one of the $D_j$'s. Therefore each component of $S \setminus {\mathcal G}(f)$ is a union of a component of $S \setminus {\mathcal Z}(f)$ with a finite number of those $D_j$'s. Define the sign of $f$ on such a component to be the one of $f$ on the corresponding component if $S \setminus {\mathcal Z}(f)$.  Now we denote the union of all components with positive (resp. negative) sign as ${C^{+}}(f)$ (resp. ${C^{-}}(f)$). As a consequence of the construction, the surfaces  ${C^{+}}(f)$ and ${C^{-}}(f)$ are incompressible. As recalled earlier, an open subset of a surface $S$ is called incompressible if the fundamental group of any of its connected components maps injectively into ${\pi_1}(S).$ The union of the connected components of ${C^{+}}(f)$ (resp. ${C^{-}}(f)$) which are neither discs nor rings is denoted by ${S^{+}}(f)$ (resp. ${S^{-}}(f)$). The surfaces ${S^{\pm}}(f)$ may be empty or
disconnected
but by construction when they are nonempty, they are incompressible.

Denote the Euler characteristic of ${S^{+}}(f)$ (resp. ${S^{-}}(f)$) by ${\chi^{+}}(f)$ (resp. ${\chi^{-}}(f)$). (we use the convention that the Euler characteristic of the empty set is zero). The incompressibility property of ${S^{+}}(f)$ and ${S^{-}}(f)$ gives that  ${\chi^{+}}(f)$ + ${\chi^{-}}(f)$ is greater than $\chi(S)$. By definition, we have ${\chi^{\pm}}(f) \leq 0$ with equality only if ${S^{\pm}}(f)$ is empty.

\textbf{Lemma 1} {\it The Euler characteristic of at least one component of $S \setminus {\mathcal G}(f)$ is negative.}

\textbf{Proof.} Let us suppose by contradiction that for some $f \neq 0 \in {\mathcal E}^{\frac{1}{4} + {\epsilon_0}(S)}$, each component $S_i$, $1 \leq i \leq m$ of $S \setminus {\mathcal G}(f)$ has non-negative Euler characteristic. So, each such component is homeomorphic either to an open disc or to an open annulus. Since $f \in {\mathcal E}^{\frac{1}{4} + {\epsilon_0}(S)}$ the Rayleigh quotient of $f$, $R(f)$ is $\leq \frac{1}{4} + {\epsilon_0}(S).$ Therefore, since $ {\mathcal G}(f)$ has measure zero, for at least one component, say $S_1$, one has $$R(f \arrowvert S_1) = \frac{\int_{S_1} {\parallel\triangledown f\parallel}^2}{\int_{S_1} {f^2}} \leq \displaystyle \frac{1}{4} + {\epsilon_0}(S).$$ Now we shall calculate the Rayleigh quotient $R(f \arrowvert S_1)$ and show that our choice of ${\epsilon_0}(S)$ leads to a contradiction.

Let us assume that $S_1$ is homeomorphic to an open disc. The case when $S_1$ is an annulus can be dealt with similarly. Since $f$ is smooth and $f|{\partial S_1} = 0$ we can choose, by Sard's theorem, a sequence $(\epsilon_n)$ of regular values of $f$ converging to $0$. Then the level set $\{x \in S_1 : f(x) = \epsilon_n\}$ is a smooth submanifold of $S$ for $n$ large enough. Furthermore one of the components of this level set confines a domain $D_n \subsetneq S_1$ homeomorphic to a closed disc with smooth boundary such that $S_1 \setminus D_n$ has arbitrarily small area. Now we consider the Rayleigh quotient, $R(f_n | D_n)$ of the function $f_n = f-{\epsilon_n}$ restricted to the region $D_n$. This function vanishes on $\partial D_n$. As $\epsilon_n$ converges to $0$, $R(f_n|D_n)$ converges to $R(f|S_1).$ Thus for any $\delta > 0$, in particular for $\displaystyle \frac{{\epsilon_0}(S)}{2}$, we can find $\epsilon_n$ small enough such that  $R(f_n|D_n) \leq \displaystyle \frac{1}{4} + {\epsilon_0}(S) + \frac{{\epsilon_0}(S)}{2} < \frac{1}{4} + 2{\epsilon_0}(S)$. Now since $D_n$ is a closed disc with smooth boundary which is contained in ${S_1} \subseteq S$, it follows from the Rayleigh quotient characterization of the first Dirichlet eigenvalue of $D_n$ that $R(f_n|D_n) \geq {\lambda_0}(D).$ By Remark \ref{rmrk} we have ${\lambda_0}(D_n) \geq  \displaystyle \frac{1}{4} + 2{\epsilon_0}(S)$. This is a contradiction when $n$ is sufficiently large. $\Box$

So some component of $S \setminus {\mathcal G}(f)$ has negative Euler characteristic. This component is a component of ${S^{\pm}}(f)$. Thus we obtain: $${\chi^{+}}(f) + {\chi^{-}}(f) < 0.$$

Now we start with some definitions and complete the proof.
\begin{definition}\label{core}
According to the sign of $f$ on $S_i$, we denote this component as ${{S_i}^{+}}(f)$ or ${{S_i}^{-}}(f)$. For each such surface with negative Euler characteristic, we consider a {\it compact core}, i.e. a compact surface ${{K_i}^\pm}(f) \subset {{S_i}^{\pm}}(f)$ such that the inclusion is a homotopy equivalence. We then define the surface ${\Sigma^{+}}(f)$ (resp. ${\Sigma^{-}}(f)$) as the union of the compact cores ${{K_i}^{+}}(f)$ (resp. ${{K_i}^{-}}(f)$) and of those components (if any) of the complement $S \setminus \bigcup {{K_i}^{+}}(f)$ (resp. $S \setminus \bigcup {{K_i}^{-}}(f)$), which are annuli. Therefore, ${\Sigma^{+}}(f)$ (resp. ${\Sigma^{-}}(f)$) is obtained from  $\bigcup {{K_i}^{+}}(f)$ (resp. $\bigcup {{K_i}^{-}}(f)$ ), by adding (if any) the annuli between the components of $\bigcup {{K_i}^{+}}(f)$ (resp. $\bigcup {{K_i}^{-}}(f)$). We call $\Sigma(f) = {\Sigma^{+}}(f) \bigcup {\Sigma^{-}}(f)$, the characteristic surface of $f$, while ${\Sigma^{+}}(f)$ (resp. ${\Sigma^{-}}(f)$) is called the positive (resp. negative) characteristic surface of $f$. The definition of these surfaces depend uniquely on the choice of compact cores and those are well defined up to isotopy. By construction the Euler characteristic of  ${\Sigma^{+}}(f)$ (resp. ${\Sigma^{-}}(f)$) is ${\chi^{+}}(f)$ (resp. ${\chi^{-}}(f)$). It is clear that  ${\Sigma^{+}}(-f) = {\Sigma^{-}}(f)$ and ${\Sigma^{-}}(-f) = {\Sigma^{+}}(f)$.
\end{definition}
\subsubsection*{Continuation of the proof of Theorem 1.}
Let $m$ denote the dimension of the space ${\mathcal E}^{\frac{1}{4} + {\epsilon_0}(S)}$. Theorem 1 will follow from the inequality $m \leq (2g-2)$. Let $\mathbb S({\mathcal E}^{\frac{1}{4} + {\epsilon_0}(S)})$ denote the unit sphere of ${\mathcal E}^{\frac{1}{4} + {\epsilon_0}(S)}$ (for some arbitrary norm) and let $\mathbb P({\mathcal E}^{\frac{1}{4} + {\epsilon_0}(S)})$ be the projective space of ${\mathcal E}^{\frac{1}{4} + {\epsilon_0}(S)}$ i.e. the quotient of $\mathbb S({\mathcal E}^{\frac{1}{4} + {\epsilon_0}(S)})$ by the involution $ f \rightarrow$ -$f$.
$\\$For each integer $i$ with ${2-2g} \leq i \leq {-1}$, we denote $$C_i = \{ f \in \mathbb S({\mathcal E}^{\frac{1}{4} + {\epsilon_0}(S)})   \rvert    {\chi^{+}}(f) + {\chi^{-}}(f) = i \}.$$ According to the lemma and its consequence above, $\mathbb S({\mathcal E}^{\frac{1}{4} + {\epsilon_0}(S)}) = {\bigcup_{2-2g}^{-1}}  C_i$. On the other hand, each $C_i$ is invariant under the antipodal involution. Let $P_i$ be the quotient of $C_i$ under the antipodal involution. The projective space $\mathbb P({\mathcal E}^{\frac{1}{4} + {\epsilon_0}(S)})$ is the union of the sets $P_i$.

\textbf{Lemma 2} {\it For any integer $i$, $2-2g \leq i \leq -1$, the covering map $ C_i \rightarrow P_i$ is trivial.}

\textbf{Proof.} Let $f\in C_i.$ We use the notations introduced in the definition of characteristic surface of $f$:  ${{S_i}^{\pm}}(f)$ is a connected component of negative Euler characteristic of ${S^{\pm}}(f)$ and ${{K_i}^{\pm}}(f)$ is a compact core of ${{S_i}^{\pm}}(f)$. We may assume that the compact core has been chosen in such a way that any connected component of $Z(f)$ that is contained in some ${{S_i}^{\pm}}(f)$ is indeed contained in the interior of the corresponding ${{K_i}^{\pm}}(f)$.

For any function $g \in {\mathcal E}^{\frac{1}{4} + {\epsilon_0}(S)}$ close enough to $f$, and for each $i$, ${{K_i}^{\pm}}(f)$ is contained in a component ${{S_l}^{\pm}}(g)$ of ${{S}^{\pm}}(g)$. Fix a neighborhood $V(f)$ of $f$ in $\mathbb S({\mathcal E}^{\frac{1}{4} + {\epsilon_0}(S)})$ such that these inclusions occur on each surface ${{K_i}^{\pm}}(f)$.

We will show that for any $g \in C_i \cap V(f)$, the characteristic surfaces ${{\Sigma}^{+}}(f)$ and ${{\Sigma}^{+}}(g)$ (resp. ${{\Sigma}^{-}}(f)$ and ${{\Sigma}^{-}}(g)$) are isotopic. Choose the compact cores ${{K_i}^{\pm}}(g)$ of surfaces ${{S}^{\pm}}(g)$ so that when ${{K_i}^{\pm}}(f)$ is contained in ${{S_j}^{\pm}}(g)$, it is also contained in the interior of ${{K_j}^{\pm}}(g)$. Now observe that if two components of the boundaries of surfaces ${{K_j}^{+}}(f)$ are homotopic in $S$ then the homotopy between them is achieved by an annulus contained in ${{\Sigma}^{+}}(f)$, by the definition of the characteristic surface. Since this annulus joins two curves of ${{K_j}^{+}}(g)$ by the definition of the characteristic surface again, it is contained in one of the connected components ${{\Sigma}^{+}}(g)$ too.

We deduce from this that each connected component of ${{\Sigma}^{\pm}}(f)$ is contained in a connected component of ${{\Sigma}^{\pm}}(g)$ (of the same sign). Since ${{\Sigma}^{+}}(f)$ and ${{\Sigma}^{-}}(f)$ are incompressible in $S$, they are incompressible in ${{\Sigma}^{+}}(g)$ and ${{\Sigma}^{-}}(g)$ respectively. In particular, their Euler characteristic satisfy $${\chi^{+}}(f) \leq {\chi^{+}}(g) ~ \textrm{and} ~ {\chi^{-}}(f)) \leq {\chi^{-}}(g);$$ these inequalities can be equalities if and only if the surfaces ${\Sigma^{+}}(f)$ and ${\Sigma^{+}}(g)$ (resp. ${\Sigma^{+}}(f)$ and ${\Sigma^{+}}(g)$) are isotopic. But since $g \in C_i$, we have $${{\chi}^{+}}(f) + {{\chi}^{-}}(f) = i = {{\chi}^{+}}(g) + {{\chi}^{-}}(g).$$ Thus ${{\Sigma}^{+}}(f)$ and ${{\Sigma}^{+}}(g)$ are isotopic. The same holds for ${{\Sigma}^{-}}(f)$ and ${{\Sigma}^{-}}(g)$.

Since the {\it isotopy class of} ${{\Sigma}^{+}}(f)$ and {\it isotopy class of} ${{\Sigma}^{-}}(f)$ are locally constant on $C_i$, they are constant on each connected component of $C_i$. Finally we observe that the functions $f$ and $-f$ can not be in the same connected component of $C_i$. This is because then ${{\Sigma}^{+}}(f)$ and ${{\Sigma}^{-}}(f)$ would be isotopic. But two disjoint and incompressible surfaces of negative Euler characteristic contained in $S$ can not be isotopic. Thus the covering map in Lemma 2 is trivial. $\Box$

\subsubsection*{Continuation of the proof of Theorem 1.}
We conclude the proof of the Theorem following a method of B. S\'{e}vennec \cite{Se}. The double covering $\mathbb S({\mathcal E}^{\frac{1}{4} + {\epsilon_0}(S)}) \rightarrow \mathbb P({\mathcal E}^{\frac{1}{4} + {\epsilon_0}(S)})$ is associated to a cohomology class $\beta \in {H^1}(\mathbb P({\mathcal E}^{\frac{1}{4} + {\epsilon_0}(S)}), {\mathbb Z}/ {2\mathbb Z} )$. Each covering $C_i \rightarrow P_i$ is described by the Cech cohomology class, $\beta \rvert_{P_i}$. Since each of this covering is trivial, we have $\beta \rvert_{P_i} = 0$. Since $\mathbb P({\mathcal E}^{\frac{1}{4} + {\epsilon_0}(S)})$ is the union of $P_i$ and since there are at most $2g-2$ of them, we have: ${\beta}^{2g-2} = 0 $(\cite[Lemma 8]{Se}). Since $\beta$ has order $m$ in the ${\mathbb Z}/ {2\mathbb Z}$-cohomology ring of $\mathbb P({\mathcal E}^{\frac{1}{4} + {\epsilon_0}(S)})$, we have $ m \leq 2g-2. \Box$

\section{Systole and the Laplace spectrum.}
In this section we study the eigenvalues of the Laplacian as functions on the moduli space. Recall that the moduli space ${\mathcal M}_g$ is the space of all closed hyperbolic surfaces of genus $g$ up to isometry. ${\mathcal M}_g$ can be compactified to a space $\overline{{\mathcal M}_g}$ by adding the moduli spaces of (not necessarily connected) non-compact finite area hyperbolic surfaces with area equal to $4 \pi(g-1)$. In this compactification a sequence $(S_n)$ in ${\mathcal M}_g$, with $s(S_n) \rightarrow 0$, converges to ${S_\infty} \in {\mathcal M}_{{g_0}, {n_0}}$ (with $2{g_0}-2+{n_0} = 2g-2$) if and only if for any given $\epsilon > 0$ the $\epsilon$-thick part $({S_n}^{[\epsilon, \infty)})$ converge to ${S_\infty}^{[\epsilon, \infty)}$ in the Gromov-Hausdorff topology. Recall that the $\epsilon$-thick part of a surface $S$ is the subset of those points of $S$ where the {\it injectivity radius} is at least $\epsilon$. Recall also that the injectivity radius of a point $p \in S$ is the radius of the largest geodesic disc that can be embedded in $S$ with center $p$.

It is a classical result that for any $i$, $\lambda_i$ is a continuous function on ${\mathcal M}_g$ (see for instance \cite{C-C}). It is also shown in \cite{C-C} that eigenvalues less than 1/4 are continuous up to $\partial {\mathcal M}_g$. Let $(S_n)$ be a sequence of surfaces in ${\mathcal M}_g$ which tends to ${S_{\infty}} \in {\partial {\mathcal M}_g}= \overline{{\mathcal M}_g} \setminus {\mathcal M}_g$.

\textbf{Theorem (\cite{C-C}, \cite{He})} {\it Let $\lambda(S_n)$ be a sequence of eigenvalues of $S_n$ which converges to $\lambda < 1/4$. Then $\lambda$ is an eigenvalue of $S_{\infty}$ and up to extracting a subsequence and possibly multiplying by a scaling constant the corresponding eigenfunctions on $S_n$ converge to an eigenfunction on $S_{\infty}$ uniformly over compact subsets.}

Our situation is a bit different. For a fixed $i$ we shall study the behavior of ${\lambda_i}(S_n)$ when ${S_n} \in {\mathcal M}_g$ tends to a point in ${\partial {\mathcal M}_g}$. Recall $$V_i = \{ {\lim_{n \rightarrow \infty}} {{\lambda_i}(S_{n})} :({S_n}) ~ \textrm{is a sequence in} ~ {{\mathcal M}_g} ~ \textrm{converging to a point in} ~ {\partial {\mathcal M}_g} $$ $$ \textrm{such that the limit exists}\}.$$ In \cite{R3} Randol showed a limiting behavior of $\lambda_{2g-2}$ over some special family. Now we apply Theorem 1 to prove the following,

\textbf{Claim 1}  {\it ${\lambda_{2g-2}}(S_n)$ tends to $\displaystyle \frac{1}{4}$ if and only if $s(S_n)$ tends to zero. In particular $V_{2g-2} = \{ \displaystyle \frac{1}{4} \}$.}

\textbf{Proof.} By Theorem 1, if ${\lambda_{2g-2}}(S_n)$ tends to $\frac{1}{4}$ then ${\epsilon_0}(S_n)$ tends to zero. For the other direction we use Buser's construction. By the  definition of the systole, there is a closed geodesic $\tau$ on $S$ such that the length of $\tau$ is equal to $s(S)$. Now from the {\it Collar Theorem} (ref. \cite{Bu}) of L. Keen\cite{K} (see also \cite{R2}) and the explicit computations in \cite[p. 219]{Bu} we see that for any $\epsilon > 0$ and any $i \geq 1$ we have $\delta >0$ such that whenever $s(S) < \delta$, we can find at least $i$ disjoint annuli in the collar neighborhood of $\tau$ of length such that the first Dirichlet eigenvalue of each of the annuli is $\leq \frac{1}{4} + \epsilon$. The corresponding eigenfunctions are orthogonal. Hence we have ${{\lambda}_{i-1}}(S) \leq \frac{1}{4} + \epsilon$. Therefore using Theorem 1 for an $i \geq 2g-1$ we obtain the convergence ${\lambda_{i}}(S_n) \rightarrow \displaystyle \frac{1}{4}. \Box$

Now we show that such a limiting behavior is not true in general for $i \leq 2g-3$. Moreover

\textbf{Claim 2} {\it For any $1 \leq i \leq 2g-3$, there exists ${\Lambda_i}(g)$, $0< {\Lambda_i}(g) \leq \displaystyle \frac{1}{4}$ such that $V_i= [0,{\Lambda_i}(g)]$.}

Before starting the proof we recall the definition of Teichm\"{u}ller space, ${\mathcal T}_g$. It is the space of all marked closed hyperbolic surfaces of genus $g$. Let $S \in {\mathcal T}_g$. Given a pair of pants decomposition of $S$, we have a coordinate system on ${\mathcal T}_g$, the Fenchel-Nielsen coordinates. ${\mathcal M}_g$ is the quotient of ${\mathcal T}_g$ by the action of ${Mod}_g$, the {\it Teichm\"{u}ller modular group}. Since ${Mod}_g$ acts properly discontinuously on ${\mathcal T}_g$, ${\mathcal T}_g$ $\rightarrow$ ${\mathcal M}_g$ is a ramified topological covering. Thus the pre-composition of this covering map with $\lambda_i$ yields a map, also denoted by $\lambda_i$, from ${\mathcal T}_g$ to $\mathbb R$. We shall use the same notation for a point in ${\mathcal T}_g$ and its image in ${\mathcal M}_g$ too.

\textbf{Proof.} We shall prove the claim for $i = 1$. The proof for $1 \leq i \leq 2g-3$ is similar. We choose a pair of pants decomposition $\mathcal P$ of a $S \in {\mathcal T}_g$ and consider the corresponding Fenchel-Nielsen coordinates $({{l^{\mathcal P}}_j}, {{\theta^{\mathcal P}}_j})$ on ${\mathcal T}_g$. Here ${{l^{\mathcal P}}_j}$'s denote the {\it length} coordinates and ${{\theta^{\mathcal P}}_j}$'s denote the {\it twist} coordinates (ref. \cite{Bu}). We fix two geodesics $\gamma$ and $\gamma ^{'}$ among the boundary geodesics of the pants decomposition $\mathcal{P}$. Thus the length functions $l_\gamma$ and $l_{\gamma^{'}}$ respectively of $\gamma$ and $\gamma ^{'}$ are among ${{l^{\mathcal P}}_j}$'s. Suppose that the pants decomposition is chosen in such a way that $\gamma$ is non-separating and $\gamma{'}$ is separating.

First we prove that $V_i$ is not empty. From a construction of P. Buser \cite[Theorem 8.1.3]{Bu} it follows that if $0 < \delta < \frac{1}{24}$ then ${\lambda_{2g-3}}(S) < \frac{1}{4}$ for any $S \in {\mathcal T}_g$ with ${{l^{\mathcal P}}_j}(S) < \delta$ for all $j$ (the number $\frac{1}{24}$ has no particular significance other than ensuring this condition). We fix one such $\delta$ and consider one $M \in {\mathcal T}_g$ such that ${{l^{\mathcal P}}_j}(M) < \delta$ for all $j$. Now consider a sequence of surfaces $({S_n}) \in {\mathcal T}_g$ such that  $({{l^{\mathcal P}}_j}, {{\theta^{\mathcal P}}_j}) (S_n) = ({{l^{\mathcal P}}_j}, {{\theta^{\mathcal P}}_j})(M)$ for all  $({{l^{\mathcal P}}_j}, {{\theta^{\mathcal P}}_j})$ except $l_\gamma$ and the ${l_\gamma}(S_n)$ coordinate decreases to zero as $n$ goes to infinity. Then $({S_n})$ converges to a point ${S_\infty} \in \partial {\mathcal M}_g$. By our choice of $\delta$ (for $M$) and since the number of components of ${S_\infty}$ is exactly one, it follows from \cite[Theorem 0.1]{C-C} that $0 < {\lim_{n \rightarrow \infty}}{\lambda_1}(S_n)= {\lambda_1}(S_\infty) = p < \frac{1}{4}$. Now consider another sequence $( {{S^{'}}_n} )$, constructed in the same way as $({S_n})$ except by varying the coordinate $l_{\gamma ^{'}}$ instead of $l_\gamma$. In this case the limiting surface of the sequence $({{S^{'}}_n})$ has two components. So using \cite{C-C} again ${\lim_{n \rightarrow \infty}}{\lambda_1}({{S^{'}}_n}) = 0$. Thus we see that $0$ and $p$ $\in V_1$, proving that $V_1$ is not empty.

Next we prove that whenever some $0< c \leq \frac{1}{4}$ is in $V_1$, the whole interval $(0, c]$ is contained in $V_1$. Since $c$ is in $V_1$ we have a sequence $({P_n})$ in ${\mathcal M}_g$ such that ${\lim_{n \rightarrow \infty}}{\lambda_1}({P_n}) = c.$ Up to extracting a subsequence, we might assume that $({P_n})$ converges to $P_{\infty} \in \partial {\mathcal M}_g$. Then $P_{\infty}$ is a finite area connected (since $c>0$) non-compact hyperbolic surface of type $({g^{'}}, m)$ (where ${g^{'}} + \frac{m}{2} = g$). For some marking of $P_n$, there is a pants decomposition of $S$, ${\gamma_1}, ..., {\gamma_k},..., {\gamma_{3g-3}}$ such that ${\gamma_1},..., {\gamma_k}$ are exactly those curves on $P_n$ whose lengths tends to zero. Consider the corresponding Fenchel-Nielsen coordinates $({l_i}, {\theta_i})_{i=1, 2,.., 3g-3}$ on ${\mathcal T}_g$. These coordinates induce coordinates on ${\mathcal T}_{{g^{'}}, m}$ which will be denoted by the same notation. In these coordinates we can choose representatives of $P_n$ in ${\mathcal T}_g$ such that $({{l_i}^n}, {\theta_i}^n)(P_n)$ converges to $({{l_i}^\infty}, {\theta_i}^\infty)$ for $i> k$ and for $i \leq k$, ${{l_i}^n}$ converges to zero. Next, using the Buser construction (\cite[Theorem 8.1.3]{Bu}), we choose a ${N_{\infty}} \in {\mathcal T}_{{g^{'}}, m}$ such that ${\lambda_1}({N_{\infty}}) = \epsilon < c.$ Then by \cite[Theorem 0.1]{C-C} for any sequence $( {N_n})$ in ${\mathcal M}_g$ converging to ${N_{\infty}}$ in $\partial {\mathcal M}_g,$ one has ${\lim_{n \rightarrow \infty}}{\lambda_1}({N_n}) = \epsilon.$ In particular we consider the sequence $(N_n)$ such that $({l_i}, {\theta_i})(N_n) = ({l_i}, {\theta_i})(N_\infty)$ for $i>k$ and $({l_i}, {\theta_i})(N_n) = ({l_i}, {\theta_i})(P_n)$ for $i\leq k$. Then ${\lim_{n \rightarrow \infty}}{\lambda_1}({N_n}) = \epsilon.$

At this point we construct a path $\sigma_n$ in ${\mathcal M}_g$ joining $P_n$ and $N_n$ for each $n$. Let us consider the path given by the coordinate axes i.e. the path first goes along the ${l_i}$ axes from ${l_i}(P_n)$  to ${l_i}(N_n)$ for each $i =k+1, k+2,..,3g-3 $ in the increasing order and then the same for $\theta_i$'s. Finally for any $t \in [\epsilon, c]$ we apply the continuity property of $\lambda_1$ on ${\mathcal M}_g$ to get a surface $Q_n$ on $\sigma_n$ such that ${\lambda_1}(Q_n)= t$. By construction each point on $\sigma$, in particular $Q_n$, has $({l_i}, {\theta_i})(Q_n) = ({l_i}, {\theta_i})(P_n)$ for $i\leq k$ and all other $({l_i}, {\theta_i})(Q_n)$ are bounded by the corresponding coordinates of $P_\infty$ and $N_\infty$. Hence $Q_n$ converges to a point ${Q_{\infty}} \in \partial {\mathcal M}_g$ and since ${\lambda_1}(Q_n)= t$ for each $n$ we have ${\lim_{n \rightarrow \infty}}{\lambda_1}({Q_n}) = t.$ Therefore $V_1$ contains $[\epsilon, c]$ and $\epsilon$ being arbitrary, the whole of $(0, c]$. That ${\Lambda_i}(g) \leq \frac{1}{4}$ follows from the last claim. $\Box$
\begin{remark} A relevant question to ask here is whether $\frac{1}{4}$ belongs to $V_k$ or not. It is not hard to show, using \cite[Theorem 0.1]{C-C}, that the existence of a surface $ N \in \partial {\mathcal M}_g$ which has $k$ components and ${\lambda_1}(N) \geq \frac{1}{4}$ would guarantee that $\frac{1}{4} \in V_k$. In this connection we would like to mention that if $g$, the genus of the surface, is large then a result of Brooks and Makover in \cite{B-M} (see also \cite{B-B-D}) shows that for any given $\epsilon$ there exist a surfaces $S \in {\mathcal M}_g$ with ${\lambda_1}(S) > \frac{3}{16} - \epsilon$ (and $> \frac{1}{4} - \epsilon$ if one assumes the Selberg conjecture).
\end{remark}
\section{Non-compact finite area hyperbolic surfaces.}
In this section we study non-compact finite area hyperbolic surfaces. Recall that ${\mathcal T }_{g,n}$ denotes the Teichm\"{u}ller space of all marked hyperbolic surfaces with finite area and of geometric type $(g,n)$. Given any pair of pants decomposition of any $S^{'} \in {\mathcal T }_{g,n}$ one can consider the Fenchel-Nielsen coordinates on ${\mathcal T }_{g,n}$. Fix one such coordinate system on ${\mathcal T }_{g,n}$. Denote by ${{\mathcal T}_{g,n}}^0$ the set of all surfaces in ${\mathcal T }_{g,n}$ all of whose twist parameters are equal to zero. Recall that each surface in ${{{\mathcal T}^0}_{g,n}}$ carries an involution $\iota$ which when restricted to each pair of pants is the orientation reversing involution that fixes the boundary components. This involution induces an involution on each eigenspace of the Laplacian. The eigenfunctions corresponding to the eigenvalue $-1$ are called {\it antisymmetric} and the corresponding eigenvalue is called an {\it antisymmetric eigenvalue}. We denote the $i$-th antisymmetric cuspidal eigenvalue of $S$ by ${{{\lambda}^{o,c}}_{i}}(S).$

We observe that in Proposition 1 we have considered domains in $S$ which are diffeomorphic either to discs or to annuli. Since $S$ is compact, the domains have compact closures. Now for $S_0 \in {\mathcal T }_{g,n}$, we may have nodal domains whose closure is not compact. To tackle this problem we consider only those domains which are diffeomorphic either to discs or to annuli and where respective boundary curves are not homotopic to puncture. For any such disc or annulus, the Cheeger's inequality is still true (ref. \cite{Cha}). The computations in Lemma \ref{annuli} then apply. Therefore for any $\Omega \subseteq S_0$ diffeomorphic either to a disc or to an annulus whose boundary curves are not homotopic to a puncture, we have an explicit constant ${\epsilon_0}(S_0) > 0 $ such that $$ {\lambda_0}(\Omega) \geq \displaystyle \frac{1}{4} + 2{\epsilon_0}(S_0) > \frac{1}{4} + {\epsilon_0}(S_0).$$

\textbf{Theorem 2} {\it For any $S_0 \in {{\mathcal T}_{g,n}}^0$ there exists an explicit constant ${\epsilon_0}(S_0) >0,$ depending only on the systole of the surface $S_0$, such that ${{{\lambda}^{o,c}}_{g}}(S_0) > \frac{1}{4} + {\epsilon_0}(S_0).$}

\textbf{Proof.} The proof proceeds along the same lines as that of Theorem 1. We choose ${\epsilon_0}(S_0)$ as above and consider ${{\mathcal E}_o}^{\frac{1}{4} + {\epsilon_0}(S_0)}$, the subspace of ${C^{\infty}}(S_0)$, spanned by the anti-symmetric cuspidal eigenfunctions with eigenvalue $\leq \frac{1}{4}  + {\epsilon_0}(S_0)$. Then we use the same arguments as in Theorem 1 to prove that the dimension of ${{\mathcal E}_o}^{\frac{1}{4} + {\epsilon_0}(S_0)}$  is less than $g$. First for $ f  \neq 0 \in  {{\mathcal E}_o}^{\frac{1}{4} + {\epsilon_0}(S_0)}$ we consider the subgraph ${\mathcal G}(f)$ of ${\mathcal Z}(f)$ obtained by suppressing those components of ${\mathcal Z}(f)$ which are bounded and homotopic to a point in $S_0$ (equivalently, those which are contained in a bounded disc in $S_0$). Next we consider the components of ${S_0} \setminus {\mathcal G}(f)$ with their signs attached as defined in \ref{core}. Denote by ${\mathcal F}(\iota)$ the fixed point set of the isometry $\iota$. The set ${\mathcal F}(\iota)$ divides $S_0$ into two isometric components ${\mathcal S}_1$ and ${\mathcal S}_2$. Each ${\mathcal S}_i$ is a non-compact finite area hyperbolic surface with geodesic boundary and genus $0$. Each puncture of $S_0$ gives rise to two ideal points , one on $\partial {\mathcal S}_1$ and another on $\partial {\mathcal S}_2$.

\textbf{Claim 3} {\it For any $f \neq 0 \in {{\mathcal E}_o}^{\frac{1}{4} + {\epsilon_0}(S_0)}$ each component of ${S_0} \setminus {\mathcal G}(f)$ is contained in one of the ${\mathcal S}_i$'s and is incompressible there.}

\textbf{Proof.} By antisymmetry of $f$ with respect to $\iota$ we have ${\mathcal F}(\iota) \subseteq {\mathcal Z}(f).$ Since each bounded component of ${\mathcal F}(\iota)$ is incompressible therefore ${\mathcal F}(\iota) \subseteq {\mathcal G}(f).$ Hence the claim follows.$\Box$

Now we can argue as in the proof of Lemma 1 to conclude that the Euler characteristic of at least one component of ${S_0} \setminus {\mathcal G}(f)$ is negative. In fact using the symmetry of ${\mathcal G}(f)$ with respect to $\iota$, the Euler Characteristic of at least one component of ${{\mathcal S}_j} \setminus {\mathcal G}(f)$ is negative for each $j =1, 2$. Next we consider the unit sphere $\mathbb S({{\mathcal E}_o}^{\frac{1}{4} + {\epsilon_0}(S_0)})$ and the projective space $\mathbb P({{\mathcal E}_o}^{\frac{1}{4} + {\epsilon_0}(S_0)})$ over ${{\mathcal E}_o}^{\frac{1}{4} + {\epsilon_0}(S_0)}$. Define ${\chi^{+}}(f)$ (respectively ${\chi^{-}}(f)$) as the sum of the Euler characteristic of the components of ${{\mathcal S}_1} \setminus {\mathcal G}(f)$ with positive sign (respectively negative). Consider the decomposition of $\mathbb S({{\mathcal E}_o}^{\frac{1}{4} + {\epsilon_0}(S_0)})$ into sets $$ {\mathcal C}_i = \{ f \in {\mathbb S}({{\mathcal E}_o}^{\frac{1}{4} + {\epsilon_0}(S_0)}):  {\chi^{+}}(f) + {\chi^{-}}(f) = i \}$$ The arguments in Lemma 2 can be applied. Using the incompressibility of components of ${{\mathcal S}_1} \setminus {\mathcal G}(f)$ the possible values of  ${\chi^{+}}(f) + {\chi^{-}}(f)$ are at most $(g-1)$ (since $\chi({\mathcal S}_i) =1-g$) for any $f \in {{\mathcal E}_o}^{\frac{1}{4} + {\epsilon_0}(S_0)}$. Exactly the same arguments as in Theorem 1 work to prove that for any integer $i$, the covering map $ C_i \rightarrow P_i$ is trivial. We conclude that the dimension of ${{\mathcal E}_o}^{\frac{1}{4} + {{\epsilon_0}(S_0)}}$ is $\leq g. \Box$


\begin{thebibliography}{}
\bibitem[B]{B} Bers, L.; Spaces of degenerating Riemann surfaces, Ann. of Math. Studies 79 (1974), 43-55.


\bibitem[B-Z]{B-Z}  Burago, Yu. D.; Zalgaller, V. A. Geometric inequalities.
Translated from the Russian by A. B. Sosinski\u{i}. Grundlehren der Mathematischen Wissenschaften, 285. Springer Series in Soviet Mathematics. Springer-Verlag, Berlin.


\bibitem[Bu]{Bu}  Buser, Peter; Geometry and spectra of compact Riemann surfaces.
Progress in Mathematics, 106. Birkh\"{a}user Boston, Inc., Boston, MA, 1992.

\bibitem[B-B-D]{B-B-D} Buser Peter, Burger Marc, Dodziuk Jozef; RIEMANN SURFACES OF LARGE GENUS AND LARGE $\lambda_1$, Geometry and Analysis on Manifolds (T. Sunada ed.), Lecture Notes in Math. 1339, Springer-Verlag, Berlin, 1988, pp. 54-63.


\bibitem[B-C]{B-C}  Buser, Peter; Courtois, Gilles; Finite parts of the spectrum of a Riemann surface. Math. Ann. 287, 523 - 530 (1990)


\bibitem[B-M]{B-M} Brooks, Robert; Makover, Eran; Riemann surfaces with large first eigenvalue. J. Anal. Math. 83 (2001), 243 - 258.

\bibitem[Cha]{Cha}   Chavel, Isaac; Eigenvalues in Riemannian geometry. Pure and Applied Mathematics, 115. Academic Press, 1984.


\bibitem[C-C]{C-C}   Colbois, Bruno; Courtois, Gilles; Les valeurs propres inf\'{e}rieures $\grave{a}$ $1/4$ des surfaces de Riemann de petit rayon d'injectivit\'{e}. Comment. Math. Helv. 64 (1989), no. 3, 349 - 362.


\bibitem[D-P-R-S]{D-P-R-S} Dodziuk J., Pignataro, Randol B., Sullivan D., Estimating small eigenvalues of Riemann surfaces, Contemp. Math. vol 64, Amer. Math. Soc. Prevedence, RI, 1987, pp 93-121.


\bibitem[H]{H}  Henrot, Antoine; Extremum problems for eigenvalues of elliptic operators. Frontiers in Mathematics. Birkh\"{a}user Verlag, Basel, 2006.


\bibitem[He]{He}  Hejhal, D. ; Regular b-groups, degenerating Riemann surfaces and spectral theory, Memoires of Amer. Math. Soc. 88, No. 437, 1990.


\bibitem[I]{I} Iwaniec, H., Introduction to the Spectral Theory of Automorphic Forms, Bibl. Rev. Mat. Iberoamericana, Revista Matem\'{a}tica Iberoamericana, Madrid, 1995.


\bibitem[K]{K} Keen, L.; Collars on Riemann surfaces, Discontinuous Groups and Riemann Surfaces, Ann. of Math. Studies No. 79, Princeton University Press, Princeton, NJ, 1974, 263 - 268.


\bibitem[O]{O}   Otal, Jean-Pierre, Three topological properties of small eigenfunctions on hyperbolic surfaces. Geometry and Dynamics of Groups and Spaces, Progr. Math. 265, Birkh\"{a}user, Bassel, 2008.


\bibitem[O-R]{O-R}   Otal, Jean-Pierre; Rosas, Eulalio; Pour toute surface hyperbolique de genre g, ${\lambda_{2g-2}}>1/4$. Duke Math. J. 150 (2009), no. 1, 101 - 115.

\bibitem[R1]{R1} Randol, B.; Small eigenvalues of the Laplace operator on compact Riemann surfaces, Bull. Amer. Math. Soc. 80 (1974), 996-1000.

\bibitem[R2]{R2} Randol, B.; Cylinders in Riemann surfaces, Comm. Math. Helv. 54, 1979, pp. 1-5.


\bibitem[R3]{R3} Randol, B.; A remark on $\lambda_{2g-2}$, Proc. Amer. Math. Soc. 108 (1990), 1081-1083.

\bibitem[Se]{Se}   S\'{e}vennec, Bruno; Multiplicity of the second Schr\"{o}dinger eigenvalue on closed surfaces. Math. Ann. 324 (2002), no. 1, 195 - 211.


\bibitem[S-W-Y]{S-W-Y}   Schoen, R.; Wolpert, S.; Yau, S. T.; Geometric bounds on the low eigenvalues of a compact surface. Geometry of the Laplace operator (Proc. Sympos. Pure Math., Univ. Hawaii, Honolulu, Hawaii, 1979), pp. 279 - 285.


\end{thebibliography}
\end{document}